\documentclass[graybox]{svmult}
\usepackage{color}
\usepackage{type1cm}        % activate if the above 3 fonts are                       % not available on your system
\usepackage{float}
\usepackage{longtable}
\usepackage{makeidx}         % allows index generation
\usepackage{graphicx}        % standard LaTeX graphics tool                             % when including figure files
\usepackage{multicol}        % used for the two-column index
\usepackage[bottom]{footmisc}% places footnotes at page bottom
\usepackage{newtxtext}       %
\usepackage{newtxmath}       % selects Times Roman as basic font
\usepackage{mathtools}
\usepackage{hyperref}
\hypersetup{
    colorlinks=true,
    linkcolor=blue,
    filecolor=magenta,
    urlcolor=cyan,
}
\usepackage[ntheorem]{empheq}
\makeindex             % used for the subject index
\numberwithin{equation}{section}
\begin{document}
\title*{BOLIB: Bilevel Optimization LIBrary\\ of test problems}
\author{Shenglong Zhou, Alain B. Zemkoho, and Andrey Tin}
\institute{Shenglong Zhou \at School of Mathematics, University of Southampton, University Road, SO17 1BJ Southampton, United Kingdom,  \email{shenglong.zhou@soton.ac.uk}
\and Alain B. Zemkoho \at School of Mathematics, University of Southampton, University Road, SO17 1BJ Southampton, United Kingdom,  \email{a.b.zemkoho@soton.ac.uk}
\and  Andrey Tin \at School of Mathematics, University of Southampton, University Road, SO17 1BJ Southampton, United Kingdom,  \email{a.tin@soton.ac.uk}}
\maketitle
%\abstract*{Each chapter should be preceded by an abstract (no more than 200 words) that summarizes the content. The abstract will appear \textit{online} at \url{www.SpringerLink.com} and be available with unrestricted access. This allows unregistered users to read the abstract as a teaser for the complete chapter.
%Please use the 'starred' version of the \texttt{abstract} command for typesetting the text of the online abstracts (cf. source file of this chapter template \texttt{abstract}) and include them with the source files of your manuscript. Use the plain \texttt{abstract} command if the abstract is also to appear in the printed version of the book.}
${}$\\[-20ex]
\abstract{This chapter presents the Bilevel Optimization LIBrary of the test
problems (BOLIB–for short), which contains a collection of test problems, with
continuous variables, to help support the development of numerical solvers for
bilevel optimization. The library contains 173 examples with 138 nonlinear, 24
linear, and 11 simple bilevel optimization problems. This BOLIB collection is
probably the largest bilevel optimization library of test problems. Moreover, as the
library is computation-enabled with the MATLAB m-files of all the examples, it
provides a uniform basis for testing and comparing algorithms. The library, together
with all the related codes, is freely available at \href{https://biopt.github.io/bolib/}{biopt.github.io/bolib}.
${}$\\[2ex]
}

\section{Introduction}\label{sec:1}
The bilevel optimization problem can take the form
\begin{equation} \label{P}
\begin{array}{rl}
   \underset{x,y}{\min}&F(x,y) \\
   \text{s.t.}&G(x,y)\leq 0,\\
              &y\in S(x):= \arg\underset{y}\min~\left\{f(x,y)\left|\; g(x,y)\leq 0\right.\right\},
\end{array}
\end{equation}
%\begin{equation}\label{P-0}
%\begin{array}{rl}
%  \underset{x, y}{\min} & F(x,y)\\
%   \text{s.t.}          & G(x,y)\leq 0,\;\; H(x,y)=0,\\[.05ex]
%                        & y\in S(x):= \arg\underset{y}{\min}\,\left\{f(x,y)\left|\;\, g(x,y)\leq 0, \;\; h(x,y)=0\right.\right\},
%\end{array}
%\end{equation}
where the functions $G :\mathbb{R}^{n_x}\times \mathbb{R}^{n_y} \rightarrow \mathbb{R}^{n_G}$ and $g :\mathbb{R}^{n_x}\times \mathbb{R}^{n_y} \rightarrow \mathbb{R}^{n_g}$ define the upper-level and lower-level constraints, respectively.
% and $g :\mathbb{R}^{n_x}\times \mathbb{R}^{n_y} \rightarrow \mathbb{R}^{n_g}$ describes the lower-level constraints.
As for $F :\mathbb{R}^{n_x}\times \mathbb{R}^{n_y} \rightarrow \mathbb{R}$ and $f :\mathbb{R}^{n_x}\times \mathbb{R}^{n_y} \rightarrow \mathbb{R}$, they  denote the upper-level and lower-level objective functions, respectively. The set-valued map $S : \mathbb{R}^{n_x} \rightrightarrows \mathbb{R}^{n_y}$ represents the optimal solution/argminimum  mapping of the lower-level problem. Further recall that problem \eqref{P} as a whole is often called upper-level problem.

Our aim is to propose a computation-enabled library of test problems to help accelerate the development of numerical methods for bilevel programs in the form \eqref{P}.
%{\color{red}{As the medium and long term goal of the library is to include various classes of bilevel optimization problems}}, in particular, simple, linear, nonlinear, and real-world applications or case studies, we intentionally consider our model \eqref{P-0} to be broad, as it is likely that the overwhelming majority of these problems will be of this form. %However,  we just focus on nonlinear academic examples here.
%%Recall that problem \eqref{P-0} is  linear if all the functions involved are linear; otherwise, it is nonlinear.
%The library is made of MATLAB-based codes of all the problems, a separate supplementary material with all the mathematical formulas, and this guideline on how to use it.
%For this version of the library, our codes are restricted to equality constraints, given that only 8\% of the problems considered have equality constraints; see Table \ref{List-example}. To deal with those problems, we use the fact that $H(x,y)=0$, for example, can be expressed as  $H(x,y)\leq0$ and $-H(x,y)\leq0$. Changes will be introduced in future versions of the library, provided that a significant number of new examples include equality constraints in the upper or lower-level of the problem. Hence, our focus here will be on bilevel optimization problems of the form
%\begin{equation} \label{P}
%\begin{array}{rl}
%   \underset{x,y}{\min}&F(x,y) \\
%   \text{s.t.}&G(x,y)\leq 0,\\
%              &y\in S(x):= \arg\underset{y}\min~\left\{f(x,y)\left|\; g(x,y)\leq 0\right.\right\}.
%\end{array}
%\end{equation}
For the sake of clarity, note that the bilevel optimization problems in this library, that we classify into the following three categories, only involve continuous variables:
\begin{itemize}
\item  \textit{Nonlinear bilevel programs}, which are problems in the form \eqref{P} with at least one of the functions involved being nonlinear.
\item  \textit{Linear bilevel programs} are problems in the form \eqref{P} with functions $F$, $f$, and all the components of $G$ and $g$ being linear.
\item \textit{Simple bilevel programs} (term coined in \cite{DempeDinhDutta2010}) are optimization problems where the feasible set is partly defined by the optimal solution set of a second optimization problem. But unlike in \eqref{P}, the lower-level problem is not a parametric optimization problem. More precisely, a simple bilevel optimization has the form:
\begin{equation} \label{SimpleP}
\begin{array}{rl}
   \underset{y}{\min}&F(y)\\
          \text{s.t.}&G(y)\leq 0,\\
                     & y\in S:= \arg\underset{y}\min~\left\{f(y)\left|\; g(y)\leq 0\right.\right\},
\end{array}
\end{equation}
where, similarly to \eqref{P}, $G :\mathbb{R}^{n_y} \rightarrow \mathbb{R}^{n_G}$  and $g :\mathbb{R}^{n_y} \rightarrow \mathbb{R}^{n_g}$ describe the upper-level and lower-level constraints, respectively,  while the real-valued function $F$ (resp. $f$), defined $\mathbb{R}^{n_y}$, represents the upper-level (resp. lower-level) objective function. The expression ``simple bilevel program'' is used in \cite{LinXuYeOnSolving2014} to refer to bilevel optimization problems of the form \eqref{P}, where $y$ (resp. $x$) is not involved in the upper-level (resp. lower-level) constraints.
%\item \textit{Parametric bilevel programs} have some functions involving extra parameters or data. We classify some examples into this category for the sake of usage for the readers.
\end{itemize}

%This paper provides a unique platform for the development of numerical methods, as well as theoretical results for bilevel optimization problems.
The main contributions of the library are three-fold. First, BOLIB provides MATLAB codes for 173 examples, including 138  nonlinear, 24 linear, and 11 simple bilevel programs, ready to be used to test numerical algorithms. Secondly, it puts together the true or best known solutions and the corresponding references for all the examples included. Hence, can serve as a benchmark platform for numerical accuracy evaluation for methods designed to solve problem \eqref{P}. Thirdly, all examples as well as their gradients and Hessians are programmed and stored in the MATLAB m-files. Thus, facilitating the use of the examples and corresponding derivatives in the implementation of numerical methods, where such information is necessary.

To the best of our knowledge, this is the largest library of test examples for bilevel optimization, especially for the nonlinear class of the problem. It includes bilevel optimization problems from Colson's BIPA \cite{C02}, Leyffer's MacMPEC \cite{MacMPEC}, as well as from Mitsos and Barton's technical report \cite{MB06}. %Special classes of examples from the latter test sets  not in this version of BOLIB will be included in future versions, where corresponding classes of problems will be expanded further.
We would like to emphasize that the fundamental objective that we hope to achieve with BOLIB is the acceleration of numerical software development for bilevel optimization, as it is our opinion that the level of expansion of applications of the problem has outpaced the development rate for numerical solvers, especially for the nonlinear class of the problem.

In the next section, we describe the library with details on the inputs and outputs of the codes, as well as some useful insights on the examples. In the subsequent section, a guideline is given on how to access the library.

\section{Description of the library}
\label{sec:2}
% Always give a unique label
% and use \ref{<label>} for cross-references
% and \cite{<label>} for bibliographic references
% use \sectionmark{}
% to alter or adjust the section heading in the running head
This section describes the structure of the library, while focusing on the inputs and outputs of each example, as well as the list of all examples together with their true or best known solutions and corresponding references. Before we proceed, note that each \texttt{m-file} contains information about the corresponding example, which include  the first and second order derivatives  of the input functions. For the upper-level objective function $F: \mathbb{R}^{n_x}\times \mathbb{R}^{n_y} \rightarrow \mathbb{R} $, these derivatives are defined as follows
%\begin{eqnarray*}
%\nabla_x F(x,y)&=&\left[
%\begin{array}{c}
%\nabla_{x_1} F\\
%\vdots\\
%\nabla_{x_{n_x}} F \\
%\end{array}\right]\in \mathbb{R}^{n_x},
%\end{eqnarray*}
\begin{equation}\label{F012}
\begin{array}{lll}
%\begin{array}{ccl}
\nabla_x F(x,y)&=&\left[
\begin{array}{c}
\nabla_{x_1} F\\
\vdots\\
\nabla_{x_{n_x}} F \\
\end{array}\right]\in \mathbb{R}^{n_x},\\[2ex]          %%%%%%%%%%%%%%%%%%%%%
\nabla^2_{xx} F(x,y)&=&\left[
\begin{array}{ccc}
\nabla^2_{x_1x_1} F& \cdots& \nabla^2_{x_{n_x}x_1} F\\
\vdots&\ddots&\vdots\\
\nabla^2_{x_1x_{n_x}} F& \cdots& \nabla^2_{x_{n_x}x_{n_x}} F\\
\end{array}\right]\in \mathbb{R}^{n_x \times n_x},\\[2ex]       %%%%%%%%%%%%%%%%%%%%%%
\nabla^2_{xy} F(x,y)&=&\left[
\begin{array}{ccc}
\nabla^2_{x_1y_1} F& \cdots& \nabla^2_{x_{n_x}y_1} F\\
\vdots&\ddots&\vdots\\
\nabla^2_{x_1y_{n_x}} F& \cdots& \nabla^2_{x_{n_x}y_{n_x}} F\\
\end{array}\right]\in \mathbb{R}^{n_y \times n_x}.
%\end{array}
\end{array}
\end{equation}
%\begin{eqnarray}\label{F012}
%%\begin{array}{ccl}
%\nabla_x F(x,y)&=&\left[
%\begin{array}{c}
%\nabla_{x_1} F\\
%\vdots\\
%\nabla_{x_{n_x}} F \\
%\end{array}\right]\in \mathbb{R}^{n_x},
%\nonumber\\
%\nabla^2_{xx} F(x,y)&=&\left[
%\begin{array}{ccc}
%\nabla^2_{x_1x_1} F& \cdots& \nabla^2_{x_{n_x}x_1} F\\
%\vdots&\ddots&\vdots\\
%\nabla^2_{x_1x_{n_x}} F& \cdots& \nabla^2_{x_{n_x}x_{n_x}} F\\
%\end{array}\right]\in \mathbb{R}^{n_x \times n_x},\\
%\nabla^2_{xy} F(x,y)&=&\left[
%\begin{array}{ccc}
%\nabla^2_{x_1y_1} F& \cdots& \nabla^2_{x_{n_x}y_1} F\\
%\vdots&\ddots&\vdots\\
%\nabla^2_{x_1y_{n_x}} F& \cdots& \nabla^2_{x_{n_x}y_{n_x}} F\\
%\end{array}\right]\in \mathbb{R}^{n_y \times n_x}.\nonumber
%%\end{array}
%\end{eqnarray}
Similar expressions are valid for $\nabla_{y} F(x,y)\in \mathbb{R}^{n_y}$, $\nabla^2_{yy} F(x,y)\in \mathbb{R}^{n_y \times n_y}$, and the lower-level objective function $f$. As the constraint functions are vector-valued, we use the following notations to refer to derivative information in the context of the upper-level constraint function $G:\mathbb{R}^{n_x}\times \mathbb{R}^{n_y} \rightarrow \mathbb{R}^{n_G}$, for instance:

\begin{eqnarray}\label{G014}
\allowdisplaybreaks
\begin{array}{ccl}
\nabla_{x} G(x,y)
&=&~
\left[
\begin{array}{c}
\nabla_{x} G_1\\
\vdots\\
\nabla_{x} G_{n_G}\\
\end{array}\right]=
\left[
\begin{array}{ccc}
\nabla_{x_1} G_{1} & \cdots& \nabla_{x_{n_x}} G_{1}\\
\vdots&\ddots&\vdots\\
\nabla_{x_1} G_{n_G} & \cdots& \nabla_{x_{n_x}} G_{n_G}
\end{array}\right] \in \mathbb{R}^{n_G \times n_x},\\[4ex]
%\end{eqnarray}
%\begin{eqnarray} \label{G013}
%\begin{array}{ccl}
\nabla^2_{xx} G(x,y)
&=&\left[
\begin{array}{c}
\nabla^2_{xx} G_1\\
\vdots\\
\nabla^2_{xx}  G_{n_G}\\
\end{array}\right]=
\left[
\begin{array}{ccc}
\nabla^2_{x_1x_1} G_1& \cdots& \nabla^2_{x_{n_x}x_1} G_1\\
\vdots&\ddots&\vdots\\
\nabla^2_{x_1x_{n_x}} G_1& \cdots& \nabla^2_{x_{n_x}x_{n_x}} G_1\\
\vdots&\ddots&\vdots\\
\nabla^2_{x_1x_1} G_{n_G} & \cdots& \nabla^2_{x_{n_x}x_1} G_{n_G}\\
\vdots&\ddots&\vdots\\
\nabla^2_{x_1x_{n_x}} G_{n_G} & \cdots& \nabla^2_{x_{n_x}x_{n_x}} G_{n_G}\\
\end{array}\right]\in \mathbb{R}^{\left(n_Gn_x\right) \times n_x},
\\[4.5ex]
\label{G014}
\nabla^2_{xy} G(x,y)
&=&\left[
\begin{array}{c}
\nabla^2_{xy} G_1\\
\vdots\\
\nabla^2_{xy}  G_{n_G}\\
\end{array}\right]=
\left[
\begin{array}{ccc}
\nabla^2_{x_1y_1} G_1& \cdots& \nabla^2_{x_{n_x}y_1} G_1\\
\vdots&\ddots&\vdots\\
\nabla^2_{x_1y_{n_y}} G_1& \cdots& \nabla^2_{x_{n_x}y_{n_y}} G_1\\
\vdots&\ddots&\vdots\\
\nabla^2_{x_1y_1} G_{n_G} & \cdots& \nabla^2_{x_{n_x}y_1} G_{n_G}\\
\vdots&\ddots&\vdots\\
\nabla^2_{x_1y_{n_y}} G_{n_G} & \cdots& \nabla^2_{x_{n_x}y_{n_y}} G_{n_G}\\
\end{array}\right]\in \mathbb{R}^{\left(n_Gn_y\right) \times n_x}.
\end{array}
\end{eqnarray}
Similar formulas are also valid for $\nabla_{y} G(x,y)\in \mathbb{R}^{n_G \times n_y}$, $\nabla^2_{yy} G(x,y)\in \mathbb{R}^{n_Gn_y \times n_y}$, and the lower-level constraint $g$. It is important to emphasize that in the context of the constraints, $\nabla_{x} G(x,y)\in \mathbb{R}^{1\times n_x}$, for example, is a row vector when $n_G=1$. However, $\nabla_{x} F(x,y)\in \mathbb{R}^{n_x}$ and $\nabla_{x} f(x,y)\in \mathbb{R}^{n_x}$ are column vectors.
%{\color{red}{It is worth noting, that for the simple examples lower-level variable is not present in the problem, hence we will only define derivatives with respect to $x$ for those examples.}}

\subsection{Inputs and outputs}
\label{subsec:2}
The \texttt{BOLIBver2} folder (see Section \ref{How to access the library} on how to access the library), which contains all the library material, includes 3 sub-folders named \texttt{Nonlinear}, \texttt{Linear}, and \texttt{Simple}. In the \texttt{Nonlinear} subfolder, there are 138 MATLAB m-files. Each one specifies a nonliner bilevel optimization test example,  named by a combination of authors' surnames, year of publication, and when necessary, the order of the example in the corresponding reference. For example, as in following figure (showing a partial list of the examples), \texttt{AiyoshiShimizu1984Ex2.m} stands for Example 2 in the paper by Aiyoshi and Shimizu published in 1984 \cite{AS84}. However, for a few examples (\texttt{DesignCentringP1}, \texttt{NetworkDesignP1}, etc.), the problem naming is based on previous use in the literature and therefore could help to easily recognize them.

\begin{center}
  \includegraphics[scale=.8]{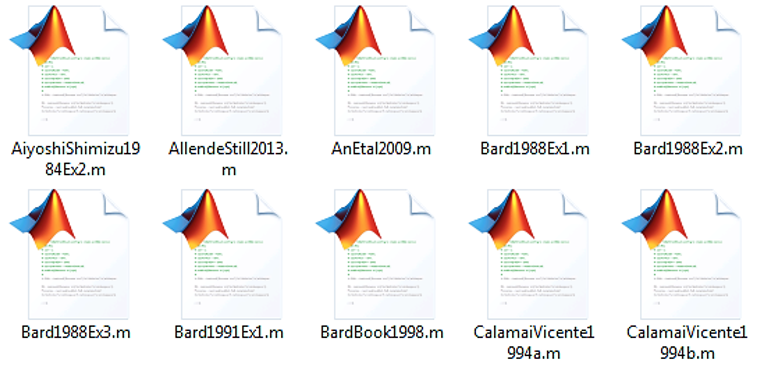}
  \end{center}
In folder \texttt{Linear}, there are 24 MATLAB m-files defining 24 liner bilevel optimization test examples. The rule of naming each example is same as in the nonlinear case. Similarly, folder \texttt{Simple} contains 11 simple bilevel optimization test examples.

 Now we describe the inputs and outputs of the m-file of a given example. Each file has the function handle named in the following way:
 \begin{eqnarray}\label{citation}
 \verb'w = example_name(x, y, keyf, keyxy)'.
 \end{eqnarray}
For the inputs, we have
\begin{eqnarray*}x&\in&\mathbb{R}^{n_x},\ \ y\in\mathbb{R}^{n_y},\\\texttt{keyf}&\in&\{ \texttt{`F', `G', `f', `g'}\},\\
\texttt{keyxy}&\in&\{[~], \texttt{`x', `y', `xx', `xy', `yy'}\},
\end{eqnarray*}
where $\texttt{`F'}$, $\texttt{`G'}$, $\texttt{`f'}$,   and $\texttt{`g'}$, respectively stand for the four functions involved in (\ref{P}). \texttt{`x'} and \texttt{`y'} represent the first order derivative with respect to $x$ and $y$, respectively. Finally, \texttt{`xx'}, \texttt{`xy'}, and \texttt{`yy'}  correspond to the second order derivative of the function $F$, $G$, $f$, and $g$, with respect to $xx$, $xy$, and $yy$, respectively.

 For the  outputs, \verb'w = example_name(x,y,keyf)' or \verb'w = example_name(x,' \verb'y,keyf,[])'  returns the function value of \texttt{keyf} while  \verb'w=example_name(x,y,keyf,' \verb'keyxy)' can additionally return the  first or second order derivative of  \texttt{keyf} w.r.t. the choice of \texttt{keyxy} as described above. We summarize all the  scenarios in Table \ref{tab:in-out}:

\begin{table}
\caption{Input--output scenarios from the m-files containing the examples}
\label{tab:in-out}
\renewcommand{\arraystretch}{1.35}\addtolength{\tabcolsep}{0pt}
\begin{tabular}{p{1.75cm}p{1.15cm}p{1.5cm}p{1.5cm}p{1.65cm}p{1.65cm}p{1.65cm}} \hline
 {\texttt{keyf}}/{\texttt{keyxy}}&$[~]$& \texttt{`x'} & \texttt{`y'} & \texttt{`xx'} &  \texttt{`xy'} &  \texttt{`yy'}\\\hline
\texttt{`F'}	&	$F(x,y)$	&	$\nabla_x F(x,y)$	&	$\nabla_y F(x,y)$	&	$\nabla^2_{xx} F(x,y)$	&	$\nabla^2_{xy} F(x,y)$	&	 $\nabla^2_{yy} F(x,y)$	\\
\texttt{`G'}	&	$G(x,y)$	&	$\nabla_x G(x,y)$	&	$\nabla_y G(x,y)$	&	$\nabla^2_{xx} G(x,y)$	&	$\nabla^2_{xy} G(x,y)$	&	 $\nabla^2_{yy} G(x,y)$	\\
\texttt{`f'}	&	$f(x,y)$	&	$\nabla_x f(x,y)$	&	$\nabla_y f(x,y)$	&	$\nabla^2_{xx} f(x,y)$	&	$\nabla^2_{xy} f(x,y)$	&	 $\nabla^2_{yy} f(x,y)$	\\
\texttt{`g'}	&	$g(x,y)$	&	$\nabla_x g(x,y)$	&	$\nabla_y g(x,y)$	&	$\nabla^2_{xx} g(x,y)$	&	$\nabla^2_{xy} g(x,y)$	&	 $\nabla^2_{yy} g(x,y)$	\\
\hline
\end{tabular}
\end{table}
For the dimension of $w$ in each scenario, see \eqref{F012}--\eqref{G014}. If $n_G=0$ (or $n_g=0$), all outputs related to $G$ (or $g$) should be empty, namely, $w=[~]$.  To further clarify the outputs, let us look at some specific usage:
\begin{itemize}
\item $w=\texttt{example\_name}(x,y,\texttt{`F'})$ or $w=\texttt{example\_name}(x,y,\texttt{`F'},{\tiny [~]})$ returns the function value of $F$, i.e., $w=F(x,y)$; this is similar for $G$, $f$, and $g$;
\item $w=\texttt{example\_name}(x,y,\texttt{`F'},\texttt{`x'})$ returns the partial derivative of $F$ with respect to $x$, i.e., $w=\nabla_x F(x,y)$;
\item $w=\texttt{example\_name}(x,y,\texttt{`G'},\texttt{`y'})$ returns the Jacobian matrix of $G$ with respect to $y$, i.e., $w=\nabla_{y} G(x,y)$;
\item $w=\texttt{example\_name}(x,y,\texttt{`f'},\texttt{`xy'})$ returns the Hessian matrix of $f$ with respect to $xy$, i.e., $w=\nabla^2_{xy} f(x,y)$;
\item $w=\texttt{example\_name}(x,y,\texttt{`g'},\texttt{`yy'})$ returns the second order derivative of $g$ with respect to $yy$, i.e., $w=\nabla^2_{yy} g(x,y)$.
\end{itemize}

 We now use two examples to illustrate the definitions above. The first one is nonlinear while the second one is a simple bilevel program.

 \begin{example}{Example} Shimizu et al. (1997), see \cite{SIB97}, considered the bilevel program \eqref{P} with
$$
\begin{array}{rll}
  F(x,y)& := & (x-5)^2 + (2y +1)^2,\\
  f(x,y)& := &(y-1)^2 -1.5xy,\\
  g(x,y)& := & \left[\begin{array}{c}
                      -3x + y+3\\
                     x-0.5y -4\\
                     x+y-7
                     \end{array}
     \right].
\end{array}
$$
\noindent Here, we have dimensions $n_x=1$, $n_y=1$, $n_G=0$, and $n_g=3$.   The \texttt{m-file} is named by \verb'ShimizuEtal1997a' (i.e., \verb'exmaple_name' = \verb'ShimizuEtal1997a' ) and was coded in MATLAB  as it can be seen in Table \ref{tab:Shimizu}. If we are given some inputs (as in left column of the table below), then \texttt{ShimizuEtal1997a} will return us corresponding results as in the right column:
$$
\begin{array}{llll}
{\rm \textbf{Inputs}} &~~~~~~~~~~&{\rm \textbf{Outputs}}\\\hline
\verb'x   ='~~ \texttt{4} &&  \verb'x   = 4'\\
\verb'y   ='~~ \texttt{0} &&  \verb'y   = 0'\\
\verb'F   ='~~ \texttt{ShimizuEtal1997a(x,y,'F')} &&  \verb'F   = 2'\\
\verb'Fx  ='~~ \texttt{ShimizuEtal1997a(x,y,'F','x')} &&  \verb'Fx  = -2'\\
\verb'Gy  ='~~ \texttt{ShimizuEtal1997a(x,y,'G','y')} && \verb'Gy  = []'\\
\verb'fxy ='~~ \texttt{ShimizuEtal1997a(x,y,'f','xy')} &&  \verb'fxy = -1.5'\\
\verb'gyy ='~~ \texttt{ShimizuEtal1997a(x,y,'g','yy')} &&  \verb'gyy = [0;0;0]'\\
\end{array}
$$
\end{example}

\begin{table}
\caption{Matlab code for \texttt{ShimizuEtal1997a}.}\label{tab:Shimizu}
{\renewcommand\baselinestretch{1.11}\selectfont
\begin{tabular}{p{0.98\textwidth}}
 \hline
 ~\\
\verb'function w=ShimizuEtal1997a(x,y,keyf,keyxy)'\\
\verb'if nargin<4 || isempty(keyxy)'\\
\verb'    switch keyf'\\
\verb'    case `F'\texttt{'}\verb'; w = (x-5)^2+(2*y+1)^2;'\\
\verb'    case `G'\texttt{'}\verb'; w = [];'\\
\verb'    case `f'\texttt{'}\verb'; w = (y-1)^2-1.5*x*y;'\\
\verb'    case `g'\texttt{'}\verb'; w = [-3*x+y+3; x-0.5*y-4; x+y-7];'\\
\verb'    end '\\
\verb'else'\\
\verb'    switch keyf'\\
\verb'    case `F'\texttt{'}\\
\verb'        switch keyxy'\\
\verb'        case `x'\texttt{'}\verb' ; w = 2*(x-5);'\\
\verb'        case `y'\texttt{'}\verb' ; w = 4*(2*y+1);'\\
\verb'        case `xx'\texttt{'}\verb'; w = 2;'\\
\verb'        case `xy'\texttt{'}\verb'; w = 0;'\\
\verb'        case `yy'\texttt{'}\verb'; w = 8;'\\
\verb'        end '\\
\verb'    case `G'\texttt{'}\\
\verb'       switch keyxy'\\
\verb'        case `x'\texttt{'}\verb' ; w = [];'\\
\verb'        case `y'\texttt{'}\verb' ; w = [];'\\
\verb'        case `xx'\texttt{'}\verb'; w = [];'\\
\verb'        case `xy'\texttt{'}\verb'; w = [];'\\
\verb'        case `yy'\texttt{'}\verb'; w = [];'\\
\verb'        end  '\\
\verb'	   case `f'\texttt{'}\\
\verb'        switch keyxy'\\
\verb'        case `x'\texttt{'}\verb' ; w = -1.5*y; '\\
\verb'        case `y'\texttt{'}\verb' ; w = 2*(y-1)-1.5*x;'\\
\verb'        case `xx'\texttt{'}\verb'; w = 0;'\\
\verb'        case `xy'\texttt{'}\verb'; w = -1.5;'\\
\verb'        case `yy'\texttt{'}\verb'; w = 2;'\\
\verb'        end           '\\
\verb'	   case `g'\texttt{'}\\
\verb'        switch keyxy'\\
\verb'        case `x'\texttt{'}\verb' ; w = [-3;  1;  1];'\\
\verb'        case `y'\texttt{'}\verb' ; w = [ 1;-0.5; 1]; '\\
\verb'        case `xx'\texttt{'}\verb'; w = [ 0;  0;  0];'\\
\verb'        case `xy'\texttt{'}\verb'; w = [ 0;  0;  0];'\\
\verb'        case `yy'\texttt{'}\verb'; w = [ 0;  0;  0];'\\
\verb'        end        '\\
\verb'     end   '\\
\verb'end'\\
\verb'end'\\
~\\
 \hline
\end{tabular} \par}
\end{table}

\newpage
\begin{example}{Example}\label{ex-FrankeEtal2018Ex513} Franke et al. (2018), see \cite{FMP18}, considered the bilevel program \eqref{P} with
$$
\begin{array}{rll}
  F(y)& := & -y_2, \\
   f(y)& := & y_3, \\
   g(y)& := &\left[\begin{array}{r}
    							y_1^2 - y_3 \\
							y_1^2 + y_2^2 - 1\\
  							 -y_3
  						\end{array}\right].\\	
\end{array}
$$
\noindent Here, we have dimensions $n_x=0$, $n_y=3$, $n_G=0$, and $n_g=3$.   The \texttt{m-file} is named by \verb'FrankeEtal2018Ex513' (i.e., \verb'exmaple_name' = \verb'FrankeEtal2018Ex513') and is equally coded in MATLAB as described in Table \ref{tab:F}.
\end{example}

\begin{table}
\caption{Matlab code for \texttt{FrankeEtal2018Ex513}. }\label{tab:F}
{\renewcommand\baselinestretch{1.11}\selectfont
\begin{tabular}{p{0.98\textwidth}}
 \hline
 ~\\
\verb'function w=FrankeEtal2018Ex513(x,y,keyf,keyxy)'\\
\verb'if nargin<4 || isempty(keyxy)'\\
\verb'    switch keyf'\\
\verb'    case `F'\texttt{'}\verb'; w = -y(2);'\\
\verb'    case `G'\texttt{'}\verb'; w = [];'\\
\verb'    case `f'\texttt{'}\verb'; w = y(3);'\\
\verb'    case `g'\texttt{'}\verb'; w = [y(1)^2-y(3); y(1)^2+y(2)^2-1; -y(3)];'\\
\verb'    end '\\
\verb'else'\\
\verb'    switch keyf'\\
\verb'    case `F'\texttt{'}\\
\verb'        switch keyxy'\\
\verb'        case `x'\texttt{'}\verb' ; w = 0;'\\
\verb'        case `y'\texttt{'}\verb' ; w = [0; -1; 0];'\\
\verb'        case `xx'\texttt{'}\verb'; w = 0;'\\
\verb'        case `xy'\texttt{'}\verb'; w = zeros(3,1);'\\
\verb'        case `yy'\texttt{'}\verb'; w = zeros(3,3);'\\
\verb'        end '\\
\verb'    case `G'\texttt{'}\\
\verb'       switch keyxy'\\
\verb'        case `x'\texttt{'}\verb' ; w = [];'\\
\verb'        case `y'\texttt{'}\verb' ; w = [];'\\
\verb'        case `xx'\texttt{'}\verb'; w = [];'\\
\verb'        case `xy'\texttt{'}\verb'; w = [];'\\
\verb'        case `yy'\texttt{'}\verb'; w = [];'\\
\verb'        end  '\\
\verb'	   case `f'\texttt{'}\\
\verb'        switch keyxy'\\
\verb'        case `x'\texttt{'}\verb' ; w = 0; '\\
\verb'        case `y'\texttt{'}\verb' ; w = [0; 0; 1];'\\
\verb'        case `xx'\texttt{'}\verb'; w = 0;'\\
\verb'        case `xy'\texttt{'}\verb'; w = zeros(3,1);'\\
\verb'        case `yy'\texttt{'}\verb'; w = zeros(3,3);'\\
\verb'        end           '\\
\verb'	   case `g'\texttt{'}\\
\verb'        switch keyxy'\\
\verb'        case `x'\texttt{'}\verb' ; w = zeros(3,1);'\\
\verb'        case `y'\texttt{'}\verb' ; w = [2*y(1) 0 -1; 2*y(1) 2*y(2) 0; 0 0 -1]; '\\
\verb'        case `xx'\texttt{'}\verb'; w = zeros(3,1);'\\
\verb'        case `xy'\texttt{'}\verb'; w = zeros(9,1);'\\
\verb'        case `yy'\texttt{'}\verb'; w = [2 0 0;0 0 0;0 0 0;2 0 0; 0 2 0;zeros(4,3)];'\\
\verb'        end        '\\
\verb'     end   '\\
\verb'end'\\
\verb'end'\\
~\\
 \hline
 \end{tabular}\par}
\end{table}

\begin{svgraybox}It is worth mentioning that despite the lack of variable $x$ in the latter example, we still treat it as an input, for the sake of unifying the inputs of the function handle as in \eqref{citation}. Hence, for all the simple bilevel optimization examples, we input \verb'x' as a scalar. In this way, \verb'x' has no impact on the example itself.
\end{svgraybox}

\subsection{Useful details on the examples}
The details related to each example presented in the BOLIB library are in a column of Table \ref{List-example} below. As we mentioned before, those examples are classified into 3 categories: nonliner, linear and simple bilevel optimisation test examples.  The first column of the table provides the list of problems, as they appear in the \texttt{Examples} subfolder of the \texttt{BOLIBver2} folder. The second column gives the reference in the literature where the example might have first appeared.   The third column combines the labels corresponding to the nature of the functions involved in \eqref{P}. Precisely, ``N'' and ``L'' will be used to indicate whether the functions $F$, $G$, $f$, and $g$ are nonlinear (N) or linear (L), while  ``$O$'' is used to symbolize that there is either no function $G$ or  $g$  present in problem \eqref{P}. Then follows the column with $n_x$ and  $n_y$ for the upper and lower-level variable dimensions, as well as $n_G$ (resp. $n_g$) to denote the  the number of components of the upper (resp. lower)-level constraint function. On the other hand, $F^*$ and $f^*$ denote the best known optimal upper and lower-level objective function values, respectively, according to the reference that is listed in the last column.

Note that examples \texttt{Zlobec2001b} and \texttt{MitsosBartonEx32} have no optimal solutions. There are 4 examples involving parameters; i.e., \texttt{CalamaiVicente1994a} with $\rho\geq 1$ (its $F^*$ and $f^*$ listed in the table are under $\rho=1$, other cases can be found in \cite{CV94}), \texttt{HenrionSurowiec2011} with $c\in\mathbb{R}$,  \texttt{IshizukaAiyoshi1992a} with $M>1$ and \texttt{RobustPortfolioP1} with $\delta\in[1,+\infty]$ (its $F^*$ and $f^*$ listed in the table are under $\delta=2$).  Dimensions $n_x$, $n_y$, $n_G$ or $n_g$ of examples \texttt{OptimalControl}, \texttt{RobustPortfolioP1}, \texttt{RobustPortfolioP2}, and \texttt{ShehuEtal2019Ex42} can be altered to get problems of different sizes, as necessary.

{\small
\begin{longtable}{p{3.5cm} p{0.6cm} |p{1.4cm}|p{0.5cm} p{0.5cm} p{0.5cm} p{0.5cm}| p{1.1cm} p{1.1cm}p{1cm} }
\caption{List of bilevel programs with related labels and known solutions.\label{List-example}}\\
  \hline
  \textbf{Example name} & \textbf{RefI} &$F$-$G$-$f$-$g$  & $n_x$ & $n_y$ & $n_G$& $n_g$&$F^*$&$f^*$&\textbf{RefII} \\ \hline
  \multicolumn{10}{c}{\textbf{Nonlinear bilevel programs}} \\ \hline
\texttt{AiyoshiShimizu1984Ex2}	&	\cite{AS84}	&	\texttt{L-L-N-L	}	&	2	&	2	&	5	&	6	&	5	&	0	&	\cite{AS84}	\\\hline
\texttt{AllendeStill2013}	&	\cite{AS13}	&	\texttt{N-L-N-N	}	&	2	&	2	&	5	&	2	&	1	&	-0.5	&	\cite{AS13}	\\\hline
\texttt{AnEtal2009}	&	\cite{ATCT09}	&	\texttt{N-L-N-L	}	&	2	&	2	&	6	&	4	&	2251.6	&	565.8	&	\cite{ATCT09}	\\\hline
\texttt{Bard1988Ex1}	&	\cite{Bard88}	&	\texttt{N-L-N-L	}	&	1	&	1	&	1	&	4	&	17	&	1	&	\cite{Bard88}	\\\hline
\texttt{Bard1988Ex2}	&	\cite{Bard88}	&	\texttt{N-L-N-L	}	&	4	&	4	&	9	&	12	&	-6600	&	54	&	\cite{C02} 	 \\\hline
\texttt{Bard1988Ex3}	&	\cite{Bard88}	&	\texttt{N-N-N-N	}	&	2	&	2	&	3	&	4	&	-12.68	&	-1.02	&	\cite{CMS05}	 \\\hline
\texttt{Bard1991Ex1}	&	\cite{Bard91}	&	\texttt{L-L-N-L	}	&	1	&	2	&	2	&	3	&	2	&	12	&	\cite{Bard91}	\\\hline
\texttt{BardBook1998Ex832}	&	\cite{BardBook98}	&	\texttt{N-L-L-L	}	&	2	&	2	&	4	&	7	&	0	&	5	&	 	 \\\hline
\texttt{CalamaiVicente1994a}	&	\cite{CV94}	&	\texttt{N-O-N-L	}	&	1	&	1	&	0	&	3	&	0&	0	&	\cite{CV94}	\\\hline
\texttt{CalamaiVicente1994b}	&	\cite{CV94}	&	\texttt{N-O-N-L	}	&	4	&	2	&	0	&	6	&	0.3125 	&	-0.4063	&	\cite{CV94}	  	 \\\hline
\texttt{CalamaiVicente1994c}	&	\cite{CV94}	&	\texttt{N-O-N-L	}	&	4	&	2	&	0	&	6	&	0.3125	&	-0.4063	&	\cite{CV94}	 \\\hline
\texttt{CalveteGale1999P1}	&	\cite{CG99}	&	\texttt{L-L-L-N	}	&	2	&	3	&	2	&	6	&	-29.2	&	0.31	&	\cite{CG99,GF01}	 \\\hline
\texttt{ClarkWesterberg1990a}	&	\cite{CW90}	&	\texttt{N-L-N-L	}	&	1	&	1	&	2	&	3	&	5	&	4	&	 \cite{PKA17}	\\\hline
\texttt{Colson2002BIPA1}	&	\cite{C02}	&	\texttt{N-L-N-L	}	&	1	&	1	&	3	&	3	&	250	&	0	&	 	 \\\hline
\texttt{Colson2002BIPA2}	&	\cite{C02}	&	\texttt{N-L-N-L	}	&	1	&	1	&	1	&	4	&	17	&	2	&	\cite{CMS05}	\\\hline
\texttt{Colson2002BIPA3}	&	\cite{C02}	&	\texttt{N-L-N-L	}	&	1	&	1	&	2	&	2	&	2	&	24.02	&	\cite{CMS05}	\\\hline
\texttt{Colson2002BIPA4}	&	\cite{C02}	&	\texttt{N-L-N-L	}	&	1	&	1	&	2	&	2	&	88.79	&	-0.77	&	\cite{CMS05}	 \\\hline
\texttt{Colson2002BIPA5}	&	\cite{C02}	&	\texttt{N-L-N-N	}	&	1	&	2	&	1	&	6	&	2.75	&	0.57	&	\cite{CMS05}	 \\\hline
%\texttt{Colson2002NDP1}	&	\cite{C02}	&	\texttt{N-L-N-L	}	&	5	&	5	&	5	&	6	&	 $\times$	&	$\times$	& 	 \\\hline
%\texttt{Colson2002NDP2}	&	\cite{C02}	&	\texttt{N-L-N-L	}	&	5	&	5	&	5	&	6	&	 $\times$	&	$\times$	& 	 \\\hline
\texttt{Dempe1992a}	&	\cite{Dempe92}	&	\texttt{L-N-N-N	}	&	2	&	2	&	1	&	2	&	$\times$	&	$\times$	&	 	 \\\hline
\texttt{Dempe1992b}	&	\cite{Dempe92}	&	\texttt{N-O-N-N	}	&	1	&	1	&	0	&	1	&	31.25	&	4	&	\cite{CMS05}	\\\hline
\texttt{DempeDutta2012Ex24}	&	\cite{DD12}	&	\texttt{N-O-N-N	}	&	1	&	1	&	0	&	1	&	0	&	0	&	\cite{DD12}	\\\hline
\texttt{DempeDutta2012Ex31}	&	\cite{DD12}	&	\texttt{L-N-N-N	}	&	2	&	2	&	4	&	2	&	-1	&	4	&	\cite{DD12}	\\\hline
\texttt{DempeEtal2012}	&	\cite{DMZ12}	&	\texttt{L-L-N-L	}	&	1	&	1	&	2	&	2	&	-1	&	-1	&	\cite{DMZ12}	\\\hline
\texttt{DempeFranke2011Ex41}	&	\cite{DF11}	&	\texttt{N-L-N-L	}	&	2	&	2	&	4	&	4	&	5	&	-2	&	\cite{DF11}	\\\hline
\texttt{DempeFranke2011Ex42}	&	\cite{DF11}	&	\texttt{N-L-N-L	}	&	2	&	2	&	4	&	3	&	2.13	&	-3.5	&	\cite{DF11}	 \\\hline
\texttt{DempeFranke2014Ex38}	&	\cite{DF14}	&	\texttt{L-L-N-L	}	&	2	&	2	&	4	&	4	&	-1	&	-4	&	\cite{DF14}	\\\hline
\texttt{DempeLohse2011Ex31a}	&	\cite{DL11}	&	\texttt{N-O-N-L	}	&	2	&	2	&	0	&	4	&	-5.5	&	0	&	\cite{DL11}	\\\hline
\texttt{DempeLohse2011Ex31b}	&	\cite{DL11}	&	\texttt{N-O-N-L	}	&	3	&	3	&	0	&	5	&	-12	&	0	& 	\\\hline
\texttt{DeSilva1978}	&	\cite{D78}	&	\texttt{N-O-N-L	}	&	2	&	2	&	0	&	4	&	-1	&	0	&	\cite{CMS05}	\\\hline
\texttt{FalkLiu1995}	&	\cite{FL95}	&	\texttt{N-O-N-L	}	&	2	&	2	&	0	&	4	&	-2.1962	&	0	&	\cite{CMS05}	 \\\hline
\texttt{FloudasEtal2013}	&	\cite{F13}	&	\texttt{L-L-N-L	}	&	2	&	2	&	4	&	7	&	0	&	200	&	\cite{TMH07}	\\\hline
\texttt{FloudasZlobec1998}	&	\cite{FZ98}	&	\texttt{N-L-L-N	}	&	1	&	2	&	2	&	6	&	1	&	-1	&	\cite{GF01,MB06}	\\\hline
\texttt{GumusFloudas2001Ex1}	&	\cite{GF01}	&	\texttt{N-L-N-L	}	&	1	&	1	&	3	&	3	&	2250	&	197.75	&	\cite{MB06}	 \\\hline
\texttt{GumusFloudas2001Ex3}	&	\cite{GF01}	&	\texttt{L-L-N-L	}	&	2	&	3	&	4	&	9	&	-29.2	&	0.31	&	\cite{MB06}	 \\\hline
\texttt{GumusFloudas2001Ex4}	&	\cite{GF01}	&	\texttt{N-L-N-L	}	&	1	&	1	&	5	&	2	&	9	&	0	&	\cite{MB06}	\\\hline
\texttt{GumusFloudas2001Ex5}	&	\cite{GF01}	&	\texttt{L-L-N-N	}	&	1	&	2	&	2	&	6	&	0.19	&	-7.23	&	\cite{MB06}	 \\\hline
\texttt{HatzEtal2013}	&	\cite{HLSB13}	&	\texttt{L-O-N-L	}	&	1	&	2	&	0	&	2	&	0	&	0	&	\cite{HLSB13}	\\\hline
\texttt{HendersonQuandt1958}	&	\cite{HQ58}	&	\texttt{N-L-N-L	}	&	1	&	1	&	2	&	1	&	-3266.7	&	-711.11	&	\cite{HQ58}	 \\\hline
\texttt{HenrionSurowiec2011}	&	\cite{HS11}	&	\texttt{N-O-N-O	}	&	1	&	1	&	0	&	0	&	$-c^2/4$	&	$-c^2/8$	&	 \cite{FJQ99}	 \\\hline
\texttt{IshizukaAiyoshi1992a}	&	\cite{IA92}	&	\texttt{N-L-L-L	}	&	1	&	2	&	1	&	5	&	0	&	-M	&	\cite{IA92}	\\\hline
\texttt{KleniatiAdjiman2014Ex3}	&	\cite{KA14}	&	\texttt{L-L-N-L	}	&	1	&	1	&	2	&	2	&	-1	&	0	&	\cite{KA14}	\\\hline
\texttt{KleniatiAdjiman2014Ex4}	&	\cite{KA14}	&	\texttt{N-N-N-N	}	&	5	&	5	&	13	&	11	&	-10	&	-3.1	&	\cite{KA14}	\\\hline
\texttt{LamparSagrat2017Ex23}	&	\cite{LS17O}	&	\texttt{L-L-N-L	}	&	1	&	2	&	2	&	2	&	-1	&	1	&	\cite{LS17O}	 \\\hline
\texttt{LamparSagrat2017Ex31}	&	\cite{LS17}	&	\texttt{N-L-L-L	}	&	1	&	1	&	1	&	1	&	1	&	0	&	\cite{LS17}	 \\\hline
\texttt{LamparSagrat2017Ex32}	&	\cite{LS17}	&	\texttt{N-O-N-O	}	&	1	&	1	&	0	&	0	&	0.5	&	0	&	\cite{LS17}	 \\\hline
\texttt{LamparSagrat2017Ex33}	&	\cite{LS17}	&	\texttt{N-L-L-L	}	&	1	&	2	&	1	&	3	&	0.5	&	0	&	\cite{LS17}	 \\\hline
\texttt{LamparSagrat2017Ex35}	&	\cite{LS17}	&	\texttt{N-L-L-L	}	&	1	&	1	&	2	&	3	&	0.8	&	-0.4	&	\cite{LS17}	 \\\hline
\texttt{LucchettiEtal1987}	&	\cite{LMP87}	&	\texttt{N-L-N-L	}	&	1	&	1	&	2	&	2	&	0	&	0	&	\cite{LMP87}	\\\hline
\texttt{LuDebSinha2016a}	&	\cite{LDS16}	&	\texttt{N-L-N-O	}	&	1	&	1	&	4	&	0	&	1.14	&	1.69	&	\cite{LDS16}	 \\\hline
\texttt{LuDebSinha2016b}	&	\cite{LDS16}	&	\texttt{N-L-N-O	}	&	1	&	1	&	4	&	0	&	0	&	1.66	&	\cite{LDS16}	 \\\hline
\texttt{LuDebSinha2016c}	&	\cite{LDS16}	&	\texttt{N-L-N-O	}	&	1	&	1	&	4	&	0	&	1.12	&	0.06	&	\cite{LDS16}	 \\\hline
\texttt{LuDebSinha2016d}	&	\cite{LDS16}	&	\texttt{L-N-L-N	}	&	2	&	2	&	11	&	3	&	$\times$	&	$\times$	&	  	 \\\hline
\texttt{LuDebSinha2016e}	&	\cite{LDS16}	&	\texttt{N-L-L-N	}	&	1	&	2	&	6	&	3	&	$\times$	&	$\times$	&	  	 \\\hline
\texttt{LuDebSinha2016f}	&	\cite{LDS16}	&	\texttt{L-N-N-O	}	&	2	&	1	&	9	&	0	&	$\times$	&	$\times$	&	   \\\hline
\texttt{MacalHurter1997}	&	\cite{MH97}	&	\texttt{N-O-N-O	}	&	1	&	1	&	0	&	0	&	81.33	&	-0.33	&	\cite{MH97}	\\\hline
\texttt{Mirrlees1999}	&	\cite{M99}	&	\texttt{N-O-N-O	}	&	1	&	1	&	0	&	0	&	1	&	0.06	&	\cite{M99}	\\\hline
\texttt{MitsosBarton2006Ex38}	&	\cite{MB06}	&	\texttt{N-L-N-L	}	&	1	&	1	&	4	&	2	&	0	&	0	&	\cite{MB06}	\\\hline
\texttt{MitsosBarton2006Ex39}	&	\cite{MB06}	&	\texttt{L-L-N-L	}	&	1	&	1	&	3	&	2	&	-1	&	-1	&	\cite{MB06}	\\\hline
\texttt{MitsosBarton2006Ex310}	&	\cite{MB06}	&	\texttt{L-L-N-L	}	&	1	&	1	&	2	&	2	&	0.5	&	-0.1	&	\cite{MB06}	\\\hline
\texttt{MitsosBarton2006Ex311}	&	\cite{MB06}	&	\texttt{L-L-N-L	}	&	1	&	1	&	2	&	2	&	-0.8	&	0	&	\cite{MB06}	\\\hline
\texttt{MitsosBarton2006Ex312}	&	\cite{MB06}	&	\texttt{N-L-N-L	}	&	1	&	1	&	2	&	2	&	0	&	0	&	\cite{MB06}	\\\hline
\texttt{MitsosBarton2006Ex313}	&	\cite{MB06}	&	\texttt{L-L-N-L	}	&	1	&	1	&	2	&	2	&	-1	&	0	&	\cite{MB06}	\\\hline
\texttt{MitsosBarton2006Ex314}	&	\cite{MB06}	&	\texttt{N-L-N-L	}	&	1	&	1	&	2	&	2	&	0.25	&	-0.08	&	\cite{MB06}	 \\\hline
\texttt{MitsosBarton2006Ex315}	&	\cite{MB06}	&	\texttt{L-L-N-L	}	&	1	&	1	&	2	&	2	&	0	&	-0.83	&	\cite{MB06}	\\\hline
\texttt{MitsosBarton2006Ex316}	&	\cite{MB06}	&	\texttt{L-L-N-L	}	&	1	&	1	&	2	&	2	&	-2	&	0	&	\cite{MB06}	\\\hline
\texttt{MitsosBarton2006Ex317}	&	\cite{MB06}	&	\texttt{N-L-N-L	}	&	1	&	1	&	2	&	2	&	0.19	&	-0.02	&	\cite{MB06}	 \\\hline
\texttt{MitsosBarton2006Ex318}	&	\cite{MB06}	&	\texttt{N-L-N-L	}	&	1	&	1	&	2	&	2	&	-0.25	&	0	&	\cite{MB06}	\\\hline
\texttt{MitsosBarton2006Ex319}	&	\cite{MB06}	&	\texttt{N-L-N-L	}	&	1	&	1	&	2	&	2	&	-0.26	&	0	&	\cite{MB06}	\\\hline
\texttt{MitsosBarton2006Ex320}	&	\cite{MB06}	&	\texttt{N-L-N-L	}	&	1	&	1	&	2	&	2	&	0.31	&	-0.08	&	\cite{MB06}	 \\\hline
\texttt{MitsosBarton2006Ex321}	&	\cite{MB06}	&	\texttt{N-L-N-L	}	&	1	&	1	&	2	&	2	&	0.21	&	-0.07	&	\cite{MB06}	 \\\hline
\texttt{MitsosBarton2006Ex322}	&	\cite{MB06}	&	\texttt{N-L-N-N	}	&	1	&	1	&	2	&	3	&	0.21	&	-0.07	&	\cite{MB06}	 \\\hline
\texttt{MitsosBarton2006Ex323}	&	\cite{MB06}	&	\texttt{N-N-L-N	}	&	1	&	1	&	3	&	3	&	0.18	&	-1	&	\cite{MB06}	\\\hline
\texttt{MitsosBarton2006Ex324}	&	\cite{MB06}	&	\texttt{N-L-N-L	}	&	1	&	1	&	2	&	2	&	-1.75	&	0	&	\cite{MB06}	\\\hline
\texttt{MitsosBarton2006Ex325}	&	\cite{MB06}	&	\texttt{N-N-N-N	}	&	2	&	3	&	6	&	9	&	-1	&	-2	&	\cite{MB06}	\\\hline
\texttt{MitsosBarton2006Ex326}	&	\cite{MB06}	&	\texttt{N-N-N-L	}	&	2	&	3	&	7	&	6	&	-2.35	&	-2	&	\cite{MB06}	\\\hline
\texttt{MitsosBarton2006Ex327}	&	\cite{MB06}	&	\texttt{N-N-N-N	}	&	5	&	5	&	13	&	13	&	2	&	-1.1	&	\cite{MB06}	\\\hline
\texttt{MitsosBarton2006Ex328}	&	\cite{MB06}	&	\texttt{N-N-N-N	}	&	5	&	5	&	13	&	13	&	-10	&	-3.1	&	\cite{MB06}	\\\hline
\texttt{MorganPatrone2006a}	&	\cite{MP06}	&	\texttt{L-L-N-L	}	&	1	&	1	&	2	&	2	&	-1	&	0	&	\cite{MP06}	\\\hline
\texttt{MorganPatrone2006b}	&	\cite{MP06}	&	\texttt{L-O-N-L	}	&	1	&	1	&	0	&	4	&	-1.25	&	0	&	\cite{MP06}	\\\hline
\texttt{MorganPatrone2006c}	&	\cite{MP06}	&	\texttt{L-O-N-L	}	&	1	&	1	&	0	&	4	&	-1	&	-0.25	&	\cite{MP06}	\\\hline
\texttt{MuuQuy2003Ex1}	&	\cite{MQ03}	&	\texttt{N-L-N-L	}	&	1	&	2	&	2	&	3	&	-2.08	&	-0.59	&	\cite{MQ03}	\\\hline
\texttt{MuuQuy2003Ex2}	&	\cite{MQ03}	&	\texttt{N-L-N-L	}	&	2	&	3	&	3	&	4	&	0.64	&	1.67	&	\cite{MQ03}	\\\hline
\texttt{NieEtal2017Ex34}	&	\cite{NWY17}	&	\texttt{L-L-N-N	}	&	1	&	2	&	2	&	2	&	2	&	0	&	\cite{NWY17}	\\\hline
\texttt{NieEtal2017Ex52}	&	\cite{NWY17}	&	\texttt{N-N-N-N	}	&	2	&	3	&	5	&	2	&	-1.71	&	-2.23	&	\cite{NWY17}	 \\\hline
\texttt{NieEtal2017Ex54}	&	\cite{NWY17}	&	\texttt{N-N-N-N	}	&	4	&	4	&	3	&	2	&	-0.44	&	-1.19	&	\cite{NWY17}	 \\\hline
\texttt{NieEtal2017Ex57}	&	\cite{NWY17}	&	\texttt{N-N-N-N	}	&	2	&	3	&	5	&	2	&	-2	&	-1	&	\cite{NWY17}	\\\hline
\texttt{NieEtal2017Ex58}	&	\cite{NWY17}	&	\texttt{N-N-N-N	}	&	4	&	4	&	3	&	2	&	-3.49	&	-0.86	&	\cite{NWY17}	 \\\hline
\texttt{NieEtal2017Ex61}	&	\cite{NWY17}	&	\texttt{N-N-N-N	}	&	2	&	2	&	5	&	1	&	-1.02	&	-1.08	&	\cite{NWY17}	 \\\hline
\texttt{Outrata1990Ex1a}	&	\cite{O90}	&	\texttt{N-O-N-L	}	&	2	&	2	&	0	&	4	&	-8.92	&	-6.05	&	\cite{O90}	\\\hline
\texttt{Outrata1990Ex1b}	&	\cite{O90}	&	\texttt{N-O-N-L	}	&	2	&	2	&	0	&	4	&	-7.56	&	-0.58	&	\cite{O90}	\\\hline
\texttt{Outrata1990Ex1c}	&	\cite{O90}	&	\texttt{N-O-N-L	}	&	2	&	2	&	0	&	4	&	-12	&	-112.71	&	\cite{O90}	\\\hline
\texttt{Outrata1990Ex1d}	&	\cite{O90}	&	\texttt{N-O-N-L	}	&	2	&	2	&	0	&	4	&	-3.6	&	-2	&	\cite{O90}	\\\hline
\texttt{Outrata1990Ex1e}	&	\cite{O90}	&	\texttt{N-O-N-L	}	&	2	&	2	&	0	&	4	&	-3.15	&	-16.29	&	\cite{O90}	\\\hline
\texttt{Outrata1990Ex2a}	&	\cite{O90}	&	\texttt{N-L-N-L	}	&	1	&	2	&	1	&	4	&	0.5	&	-14.53	&	\cite{O90}	\\\hline
\texttt{Outrata1990Ex2b}	&	\cite{O90}	&	\texttt{N-L-N-L	}	&	1	&	2	&	1	&	4	&	0.5	&	-4.5	&	\cite{O90}	\\\hline
\texttt{Outrata1990Ex2c}	&	\cite{O90}	&	\texttt{N-L-N-L	}	&	1	&	2	&	1	&	4	&	1.86	&	-10.93	&	\cite{O90}	\\\hline
\texttt{Outrata1990Ex2d}	&	\cite{O90}	&	\texttt{N-L-N-N	}	&	1	&	2	&	1	&	4	&	0.92	&	-19.47	&	\cite{O90}	\\\hline
\texttt{Outrata1990Ex2e}	&	\cite{O90}	&	\texttt{N-L-N-N	}	&	1	&	2	&	1	&	4	&	0.90	&	-14.94	&	\cite{O90}	\\\hline
\texttt{Outrata1993Ex31}	&	\cite{O93}	&	\texttt{N-L-N-N	}	&	1	&	2	&	1	&	4	&	1.56	&	-11.67	&	\cite{O93}	\\\hline
\texttt{Outrata1993Ex32}	&	\cite{O93}	&	\texttt{N-L-N-N	}	&	1	&	2	&	1	&	4	&	3.21	&	-20.53	&	\cite{O93}	\\\hline
\texttt{Outrata1994Ex31}	&	\cite{O94}	&	\texttt{N-L-N-N	}	&	1	&	2	&	2	&	4	&	3.21	&	-20.53	&	\cite{O94}	\\\hline
\texttt{OutrataCervinka2009}	&	\cite{OC09}	&	\texttt{L-L-N-L	}	&	2	&	2	&	1	&	3	&	0	&	0	&	\cite{OC09}	\\\hline
\texttt{PaulaviciusEtal2017a}	&	\cite{PKA17}	&	\texttt{N-L-N-L	}	&	1	&	1	&	4	&	2	&	0.25	&	0	&	\cite{PKA17}	 \\\hline
\texttt{PaulaviciusEtal2017b}	&	\cite{PKA17}	&	\texttt{L-L-N-L	}	&	1	&	1	&	4	&	2	&	-2	&	-1.5	&	\cite{PKA17}	 \\\hline
\texttt{SahinCiric1998Ex2}	&	\cite{SC98}	&	\texttt{N-L-N-L	}	&	1	&	1	&	2	&	3	&	5	&	4	&	\cite{SC98}	\\\hline
\texttt{ShimizuAiyoshi1981Ex1}	&	\cite{SA81}	&	\texttt{N-L-N-L	}	&	1	&	1	&	3	&	3	&	100	&	0	&	\cite{SA81}	\\\hline
\texttt{ShimizuAiyoshi1981Ex2}	&	\cite{SA81}	&	\texttt{N-L-N-L	}	&	2	&	2	&	3	&	4	&	225	&	100	&	\cite{SA81}	\\\hline
\texttt{ShimizuEtal1997a}	&	\cite{SIB97}	&	\texttt{N-O-N-L	}	&	1	&	1	&	0	&	3	&	$\times$	&	$\times$	&	 	 \\\hline
\texttt{ShimizuEtal1997b}	&	\cite{SIB97}	&	\texttt{N-L-N-L	}	&	1	&	1	&	2	&	2	&	2250	&	197.75	&	\cite{SIB97}	 \\\hline
\texttt{SinhaMaloDeb2014TP3}	&	\cite{SMD14}	&	\texttt{N-N-N-N	}	&	2	&	2	&	3	&	4	&	-18.68	&	-1.02	&	\cite{SMD14}	 \\\hline
\texttt{SinhaMaloDeb2014TP6}	&	\cite{SMD14}	&	\texttt{N-L-N-L	}	&	1	&	2	&	1	&	6	&	-1.21	&	7.62	&	\cite{SMD14}	 \\\hline
\texttt{SinhaMaloDeb2014TP7}	&	\cite{SMD14}	&	\texttt{N-N-N-L	}	&	2	&	2	&	4	&	4	&	-1.96	&	1.96	&	\cite{SMD14}	 \\\hline
\texttt{SinhaMaloDeb2014TP8}	&	\cite{SMD14}	&	\texttt{N-L-N-L	}	&	2	&	2	&	5	&	6	&	0	&	100	&	\cite{SMD14}	 \\\hline
\texttt{SinhaMaloDeb2014TP9}	&	\cite{SMD14}	&	\texttt{N-O-N-L	}	&	10	&	10	&	0	&	20	&	0	&	1	&	\cite{SMD14}	 \\\hline
\texttt{SinhaMaloDeb2014TP10}	&	\cite{SMD14}	&	\texttt{N-O-N-L	}	&	10	&	10	&	0	&	20	&	0	&	1	&	\cite{SMD14}	 \\\hline
\texttt{TuyEtal2007}	&	\cite{TMH07}	&	\texttt{N-L-L-L	}	&	1	&	1	&	2	&	3	&	22.5	&	-1.52	&	\cite{TMH07}	 \\\hline
\texttt{Vogel2002}	&	\cite{V02}	&	\texttt{N-L-N-L	}	&	1	&	1	&	2	&	1	&	1	&	-2	&	\cite{V02}	\\\hline
\texttt{WanWangLv2011}	&	\cite{ZWL11}	&	\texttt{N-O-L-L	}	&	2	&	3	&	0	&	8	&	10.63	&	-0.5	&	\cite{ZWL11}	 \\\hline
\texttt{YeZhu2010Ex42}	&	\cite{YZ10}	&	\texttt{N-L-N-L	}	&	1	&	1	&	2	&	1	&	1	&	-2	&	\cite{YZ10}	\\\hline
\texttt{YeZhu2010Ex43}	&	\cite{YZ10}	&	\texttt{N-L-N-L	}	&	1	&	1	&	2	&	1	&	1.25	&	-2	&	\cite{YZ10}	\\\hline
\texttt{Yezza1996Ex31}	&	\cite{Y96}	&	\texttt{N-L-N-L	}	&	1	&	1	&	2	&	2	&	1.5	&	-2.5	&	\cite{Y96}	\\\hline
\texttt{Yezza1996Ex41}	&	\cite{Y96}	&	\texttt{N-O-N-L	}	&	1	&	2	&	0	&	2	&	0.5	&	2.5	&	\cite{Y96}	\\\hline
\texttt{Zlobec2001a}	&	\cite{Z01}	&	\texttt{N-O-L-L	}	&	1	&	2	&	0	&	3	&	-1	&	-1	&	\cite{Z01}	\\\hline
\texttt{Zlobec2001b}	&	\cite{Z01}	&	\texttt{L-L-L-N	}	&	1	&	1	&	2	&	4	&	no	&	solution	&	\cite{Z01}	 \\\hline
\texttt{DesignCentringP1}&	 \cite{SteinStill03}	&	\texttt{N-N-N-N	}	&	3	&	6	&	3	&	3	&	 $\times$	&	$\times$	& 	 \\\hline
\texttt{DesignCentringP2}&	 \cite{SteinStill03}	&	\texttt{N-N-N-N	}	&	4	&	6	&	5	&	3	&	 $\times$	&	$\times$	& 	 \\\hline
\texttt{DesignCentringP3}&	 \cite{SteinStill03}	&	\texttt{N-N-N-N	}	&	6	&	6	&	3	&	3	&	 $\times$	&	$\times$	& 	 \\\hline
\texttt{DesignCentringP4}&	 \cite{SteinStill03}	&	\texttt{N-N-N-N	}	&	4	&	6	&	3	&	12	&	 $\times$	&	$\times$	& 	 \\\hline
\texttt{NetworkDesignP1}	&	\cite{C02}	&	\texttt{N-L-N-L	}	&	5	&	5	&	5	&	6	&	 300.5	&	419.8	& 	\cite{CMS05} \\\hline
\texttt{NetworkDesignP2}	&	\cite{C02}	&	\texttt{N-L-N-L	}	&	5	&	5	&	5	&	6	&	 142.9	&	81.95	& \cite{CMS05}	 \\\hline
\texttt{OptimalControl}	&	 \cite{MehlitzGerd16}&	\texttt{N-N-N-L	}	&	2	&	$n_y$	&	3	&	$2n_y$	&	 $\times$	&	$\times$	 	\\\hline
\texttt{RobustPortfolioP1}	&	\cite{SteinStill03}&	\texttt{L-N-N-N	}	&	N+1	&	N	&	N+3	&	N+1	&	 1.15	&	0	&\cite{SteinStill03}	\\\hline
\texttt{RobustPortfolioP2}	&	\cite{SteinStill03}&	\texttt{L-N-N-N	}	&	N+1	&	N	&	N+3	&	N+1	&	 1.15	&	0	&\cite{SteinStill03}	 	\\\hline
\texttt{TollSettingP1}	&	\cite{C02}&	\texttt{N-L-N-L	}	&	3	&	8	&	3	&	18	&	 -7	&	12	&\cite{CMS05} 	\\\hline
\texttt{TollSettingP2}	&	\cite{C02}&	\texttt{N-L-N-L	}	&	3	&	18	&	3	&	38	&	 -4.5	&	32	&\cite{CMS05} 	\\\hline
\texttt{TollSettingP3}	&	\cite{C02}&	\texttt{N-L-N-L	}	&	3	&	18	&	3	&	38	&	- 3.5	&	32	&\cite{CMS05} 	\\\hline
\texttt{TollSettingP4}	&	\cite{C02}&	\texttt{N-O-N-L	}	&	2	&	4	&	0	&	 8&	-4	&	14	&\cite{CMS05} 	\\\hline
\texttt{TollSettingP5}	&	\cite{C02}&	\texttt{N-O-N-L	}	&	1	&	4	&	0	&	8	&	-2.5	&	14	&\cite{CMS05} 	\\\hline  \multicolumn{10}{c}{\textbf{Linear bilevel programs}} \\ \hline
 \texttt{AnandalinghamWhite1990}	&	\cite{AW90}	&	\texttt{L-L-L-L	}	&	1	&	1	&	1	&	6	&	-49	&	15	&	\cite{AW90}	\\\hline
\texttt{Bard1984a}	&	\cite{Bard84}	&	\texttt{L-L-L-L	}	&	1	&	1	&	1	&	5	&	28/9	&	-60/9 	&	\cite{Bard84}	\\\hline
\texttt{Bard1984b}	&	\cite{Bard84}	&	\texttt{L-L-L-L	}	&	1	&	1	&	1	&	5	&	-37.6 	&	1.6	&	\cite{Bard84}	\\\hline
\texttt{Bard1991Ex2}	&	\cite{Bard91}	&	\texttt{L-L-L-L	}	&	1	&	2	&	1	&	5	&	-1	&	-1	&	\cite{Bard91}	\\\hline
%\texttt{BardFalk82a}	&	\cite{BF82}	&	\texttt{L-L-L-L	}	&	2	&	3	&	2	&	6	&	-26	&	3.2	&	\cite{BF82} 	 \\\hline
\texttt{BardFalk1982Ex2}	&	\cite{BF82}	&	\texttt{L-L-L-L	}	&	2	&	2	&	2	&	5	&	-3.25 	&	-4	&	\cite{BF82}	 \\\hline
\texttt{Ben-AyedBlair1990a}	&	\cite{BAB90}	&	\texttt{L-L-L-L	}	&	1	&	2	&	2	&	4	&	-2.5	&	-5	&	\cite{BAB90}	\\\hline
\texttt{Ben-AyedBlair1990b}	&	\cite{BAB90}	&	\texttt{L-L-L-L	}	&	1	&	1	&	1	&	4	&	-6	&	5	&	 \cite{BAB90}	 \\\hline
\texttt{BialasKarwan1984a}	&	\cite{BK84}	&	\texttt{L-L-L-L	}	&	1	&	2	&	1	&	7	&	-2 &	 -0.5	&	\cite{BK84}	\\\hline
\texttt{BialasKarwan1984b}	&	\cite{BK84}	&	\texttt{L-L-L-L	}	&	1	&	1	&	1	&	6	&	-11 	&	11	&	\cite{BK84}	  	 \\\hline
\texttt{CandlerTownsley1982}	&	\cite{CT82}	&	\texttt{L-L-L-L	}	&	2	&	3	&	2	&	6	&	-29.2	 &	3.2	&	\cite{CT82}	 \\\hline
\texttt{ClarkWesterberg1988}	&	\cite{CW88}	&	\texttt{L-O-L-L	}	&	1	&	1	&	0	&	3	&	-37	&	14	&	\cite{CW88}	 \\\hline
\texttt{ClarkWesterberg1990b}	&	\cite{CW90}	&	\texttt{L-L-L-L	}	&	1	&	2	&	2	&	5	&	-13	&	-4	&	 \cite{CW90}	\\\hline
\texttt{GlackinEtal2009}	&	\cite{GlackinEckerKupferschmid2009}	&	\texttt{L-L-L-L	}	&	2	&	1	&	3	&	3	&	6	&	0	&	 \cite{GlackinEckerKupferschmid2009}	 \\\hline
\texttt{HaurieSavardWhite1990}	&	\cite{HSW90}	&	\texttt{L-O-L-L	}	&	1	&	1	&	0	&	4	&	27	&	-3	&	\cite{HSW90}	\\\hline
\texttt{HuHuangZhang2009}	&	\cite{HHZ2009}	&	\texttt{L-L-L-L	}	&	1	&	2	&	1	&	5	&	-76/9	 &	-41/9 	&	\cite{HHZ2009}	\\\hline
\texttt{LanWenShihLee2007}	&	\cite{LanWenShihLee07}	&	\texttt{L-L-L-L	}	&	1	&	1	&	1	&	7	&	-85.09	&	50.17	&	\cite{LanWenShihLee07}	 \\\hline
\texttt{LiuHart1994}	&	\cite{LH94}	&	\texttt{L-L-L-L	}	&	1	&	1	&	1	&	4	&	-16	&	4	&	\cite{LH94}	\\\hline
\texttt{MershaDempe2006Ex1}	&	\cite{MD06}	&	\texttt{L-L-L-L	}	&	1	&	1	&	1	&	5	&	$\times$	&	$\times$	&	\cite{MD06}	\\\hline
\texttt{MershaDempe2006Ex2}	&	\cite{MD06}	&	\texttt{L-L-L-L	}	&	1	&	1	&	2	&	2	&	-20	&	-6	&	\cite{MD06}	\\\hline
\texttt{TuyEtal1993}	&	\cite{Tuy1993}	&	\texttt{L-L-L-L	}	&	2	&	2	&	3	&	4	&	-3.25 	&	-6	&	\cite{Tuy1993}	 \\\hline
\texttt{TuyEtal1994}	&	\cite{Tuy94}	&	\texttt{L-L-L-L	}	&	2	&	2	&	3	&	3	&	6 	&	0	&	\cite{Tuy94}	 \\\hline
\texttt{TuyEtal2007Ex3}	&	\cite{TMH07}	&	\texttt{L-L-L-L	}	&	10	&	6	&	12	&	13	&	-467.46 	&	-11.62	&	\cite{TMH07}	 \\\hline
\texttt{VisweswaranEtal1996}	&	\cite{VFI96}	&	\texttt{L-L-L-L	}	&	1	&	1	&	1	&	5	&	28/9	&	-60/9	&	 \cite{VFI96}	 \\\hline
\texttt{WangJiaoLi2005}	&	\cite{WangJiaoLi2005}	&	\texttt{L-L-L-L	}	&	1	&	2	&	2	&	2	&	-1000	&	-1	&	 \cite{WangJiaoLi2005}	 \\\hline
  \multicolumn{10}{c}{\textbf{Simple bilevel programs}} \\ \hline
\texttt{FrankeEtal2018Ex53}	&	\cite{FMP18}	&	\texttt{N-L-N-L	}	&	0	&	2	&	4	&	4	&	1	&	1 	&	\cite{FMP18}	\\\hline
\texttt{FrankeEtal2018Ex511}	&	\cite{FMP18}	&	\texttt{N-O-L-L	}	&	0	&	3	&	0	&	4	&	3	&	0 	&	\cite{FMP18}	\\\hline
\texttt{FrankeEtal2018Ex513}	&	\cite{FMP18}	&	\texttt{L-O-L-N	}	&	0	&	3	&	0	&	3	&	-1	&	0 	&	\cite{FMP18}	\\\hline
\texttt{FrankeEtal2018Ex521}	&	\cite{FMP18}	&	\texttt{L-O-L-N	}	&	0	&	2	&	0	&	3	&	-1	&	0 	&	\cite{FMP18}	\\\hline
  \texttt{MitsosBarton2006Ex31}	&	\cite{MB06}	&	\texttt{L-L-L-L	}	&	0	&	1	&	2	&	2	&	1	&	-1 	&	\cite{MB06}	\\\hline
\texttt{MitsosBarton2006Ex32}	&	\cite{MB06}	&	\texttt{L-L-L-L	}	&	0	&	1	&	3	&	2	&	no	&	solution 	&	\cite{MB06}	\\\hline
\texttt{MitsosBarton2006Ex33}	&	\cite{MB06}	&	\texttt{L-L-N-N	}	&	0	&	1	&	2	&	3	&	-1	&	1 	&	\cite{MB06}	\\\hline
\texttt{MitsosBarton2006Ex34}	&	\cite{MB06}	&	\texttt{L-L-N-L	}	&	0	&	1	&	2	&	2	&	1	&	-1 	&	\cite{MB06}	\\\hline
\texttt{MitsosBarton2006Ex35}	&	\cite{MB06}	&	\texttt{L-L-N-L	}	&	0	&	1	&	2	&	2	&	0.5	&	-1 	&	\cite{MB06}	\\\hline
\texttt{MitsosBarton2006Ex36}	&	\cite{MB06}	&	\texttt{L-L-N-L	}	&	0	&	1	&	2	&	2	&	-1	&	-1 	&	\cite{MB06}	\\\hline
\texttt{ShehuEtal2019Ex42}	&	\cite{SVZ19}	&	\texttt{N-O-N-O	}	&	0 & $n_y$ & 0& 0	&	$\times$ 	&	 $\times$ 	&	 	\\\hline
%  \multicolumn{10}{c}{\textbf{Parametric bilevel examples}} \\ \hline
%  \texttt{CalamaiVicente1994a}	&	\cite{CV94}	&	\texttt{N-O-N-L	}	&	1	&	1	&	0	&	3	&	 &	 	&	\cite{CV94}	\\\hline
%\texttt{HenrionSurowiec2011}	&	\cite{HS11}	&	\texttt{N-O-N-O	}	&	1	&	1	&	0	&	0	&	$-c^2/4$	&	$-c^2/8$	&	 \cite{FJQ99}	 \\\hline	
%\texttt{IshizukaAiyoshi1992a}	&	\cite{IA92}	&	\texttt{N-L-L-L	}	&	1	&	2	&	1	&	5	&	0	&	-$M$	&	\cite{IA92}	\\\hline
%\texttt{IOC1992}	&	 \cite{MehlitzGerd16}&	\texttt{N-N-N-L	}	&	2	&	$n_y$	&	3	&	$2n_y$	&	 &	 	\\\hline
%\texttt{ShehuEtal2019Ex42}	&	\cite{SVZ19}	&	\texttt{N-O-N-O	}	&	$n_x$ & $n_y$ & 0& 0	&	 	&	  	&	 	\\\hline
  \end{longtable}}

\begin{svgraybox}  It is worth pointing out that some examples involve equalities constraints in the upper or lower-level problems. As only 8\% of the BOLIB problems have such constraints, we preserve the uniformity in the structure of the codes by converting equalities constraints into inequalities. For the sake of clarity, we list all the examples with equality constraints below.
\end{svgraybox}

{\small
\begin{longtable}{p{3.5cm} p{1.45cm} |p{2.0cm}|p{0.6cm} p{0.6cm} p{0.6cm}p{0.6cm}p{0.6cm} p{0.6cm}}
\caption{List of  bilevel programs with equalities constraints.\label{List-example-equ}}\\
  \hline
  \textbf{Example name} & \textbf{Ref.} &$F$-$G$-$H$-$f$-$g$-$h$  & $n_x$ & $n_y$ & $n_G$& $n_H$& $n_g$& $n_h$ \\ \hline

  \texttt{DempeDutta2012Ex31} & \cite{DD12} & \texttt{L-L-N-N-N-O } & 2 & 2 & 2 &  1&2& 0  \\\hline
  \texttt{DempeFranke2011Ex41} & \cite{DF11} & \texttt{N-L-L-N-L-O } & 2 & 2 & 2 &  1&4& 0  \\
  \hline
  \texttt{DempeFranke2011Ex42} & \cite{DF11} & \texttt{N-L-L-N-L-O } & 2 & 2 & 2 &  1 &3& 0 \\
  \hline

  \texttt{Zlobec2001b} & \cite{Z01} &\texttt{L-L-O-L-L-N} & 1 & 1 & 2 &  0& 2 & 1   \\
  \hline
  \texttt{NetworkDesignP1}	&	\cite{C02}	 & \texttt{N-L-O-N-O-L } & 5 & 5 & 5 &  0&0& 3  \\\hline
  \texttt{NetworkDesignP2}	&	\cite{C02}	 & \texttt{N-L-O-N-O-L } & 5 & 5 & 5 &  0&0& 3  \\\hline
    \texttt{OptimalControl}	&	 \cite{MehlitzGerd16}&	\texttt{N-N-O-N-L-L}	&	2	&	$n_y$	&	3	&0&	$n_y$	&	$\frac{1}{2}n_y$  	\\\hline
\texttt{RobustPortfolioP1}	&	\cite{SteinStill03}&	\texttt{L-N-L-N-N-O	}	&	N+1	&	N	&	N+1	& 1&	N+1	&	 0	 	\\\hline
\texttt{RobustPortfolioP2}	&	\cite{SteinStill03}&	\texttt{L-N-L-N-N-O	}	&	N+1	&	N	&	N+1	& 1&	N+1	&	 0	 	\\\hline
 \texttt{TollSettingP1}	&	\cite{C02} &	\texttt{N-L-O-N-L-L	}	&	3	&	8	&	3& 0 &	8	&	5	 	\\\hline
  \texttt{TollSettingP2}	&	\cite{C02} &	\texttt{N-L-O-N-L-L	}	&	3	&	18	&	3& 0 &	18	&	10	 	\\\hline
    \texttt{TollSettingP3}	&	\cite{C02} &	\texttt{N-L-O-N-L-L	}	&	3	&	18	&	3& 0 &	18	&	10	 	\\\hline
    \texttt{TollSettingP4}	&	\cite{C02} &	\texttt{N-O-O-N-L-L	}	&	2	&	 4	&	0& 0 &	2	&	4	 	\\\hline
    \texttt{TollSettingP5}	&	\cite{C02} &	\texttt{N-O-O-N-L-L	}	&	1	&	 4	&	0& 0 &	2	&	4	 	\\\hline
    \end{longtable}}

\section{How to access the library?}\label{How to access the library}
%The BOLIB library (versions 1 and 2) can be accessed through the dedicated website \href{https://biopt.github.io/bolib/}{biopt.github.io/bolib}. Under this link, you will find the zipped folder named \texttt{BOLIBver2} containing all the relevant files for the version of the library presented in this paper. The folder contains the subfolder named \texttt{Examples}, which contains all the m-files with the codes of the examples, as described in the previous section. The pdf file named \texttt{Formulas} collects all the mathematical formulas of all the examples in this library. To start with the library, it is advised to consult the \texttt{readme} file for some useful instructions on how to use it.
The library can be accessed through the dedicated website \href{https://biopt.github.io/bolib/}{biopt.github.io/bolib}. Under this link, you will find the zipped folder named BOLIBver2
containing all the relevant files for the version of the library presented in this paper.
The folder contains the subfolder named \texttt{Examples}, which contains all the m-files
with the codes of the examples, as described in the previous section. The pdf file
named \texttt{Formulas} collects all the mathematical formulas of all the examples in
this library. To start with the library, it is advised to consult the \texttt{readme} file for
some further useful instructions on how to use it.

\begin{acknowledgement}
The work of the first and second authors is partly funded by the  EPSRC Grant {EP/P022553/1}. The third author's work is partly funded by the University of Southampton's Presidential Scholarship. We thank  Dr  Patrick Mehlitz (Brandenburgische Technische Universit\"{a}t Cottbus-Senftenberg) for the \texttt{OptimalControl} example and related MATLAB files.
\end{acknowledgement}
%%
%\section*{Appendix}
%\addcontentsline{toc}{section}{Appendix}
%%
%%
%When placed at the end of a chapter or contribution (as opposed to at the end of the book), the numbering of tables, figures, and equations in the appendix section continues on from that in the main text. Hence please \textit{do not} use the \verb|appendix| command when writing an appendix at the end of your chapter or contribution. If there is only one the appendix is designated ``Appendix'', or ``Appendix 1'', or ``Appendix 2'', etc. if there is more than one.
%
%\begin{equation}
%a \times b = c
%\end{equation}

\end{document}